\documentclass[12pt]{amsart}
\usepackage{amsfonts, amsmath, latexsym, epsfig}
\usepackage{amssymb}
\usepackage{epsf}
\usepackage{url}

\newcommand{\sfA}{\ensuremath{\mathsf{A}}}
\newcommand{\RR}{\ensuremath{\mathbb{R}}}

\newcommand{\ZZ}{\ensuremath{\mathbb{Z}}}

\newcommand{\N}{\ensuremath{\mathbb{N}}}

\newtheorem{theorem}{Theorem}

\newtheorem{definition}{Definition}
\def\QuotS#1#2{\leavevmode\kern-.0em\raise.2ex\hbox{$#1$}\kern-.1em/\kern-.1em\lower.25ex\hbox{$#2$}}

\DeclareMathOperator{\vertt}{vert}
\DeclareMathOperator{\conv}{conv}
\DeclareMathOperator{\SC}{SC}

\begin{document}

\author{Mathieu Dutour Sikiri\'c}
\address{Mathieu Dutour Sikiri\'c, Rudjer Boskovi\'c Institute, Bijenicka 54, 10000 Zagreb, Croatia}
\email{mdsikir@irb.hr}

\author{Konstantin Rybnikov}
\address{K. Rybnikov, Department of Mathematical Sciences, University of Massachusetts at Lowell, Lowell, MA 01854, USA}
\email{Konstantin\_Rybnikov@uml.edu}

\thanks{The authors are thankful for the hospitality of the Hausdorff Research Institute for Mathematics in Bonn, where this research was done. First author has been supported by the Croatian Ministry of Science, Education and Sport under contract 098-0982705-2707.}

\title{Perfect but not generating Delaunay polytopes}
\date{}

\maketitle

\begin{abstract} In his seminal 1951 paper ``Extreme forms" Coxeter \cite{cox51}
observed that for $n \ge 9$ one can add vectors to the
perfect lattice $\sfA_9$ so that the resulting perfect lattice, called
$\sfA_9^2$ by Coxeter, has exactly the same set of minimal vectors.
An inhomogeneous analog of the notion of perfect lattice is that of
a lattice with a perfect Delaunay polytope: 
the vertices of a perfect Delaunay polytope are
the analogs of minimal vectors in a perfect lattice.
We find a new infinite series $P(n,s)$ for $s\geq 2$
and $n+1\geq 4s$ of
$n$-dimensional perfect Delaunay polytopes.
A remarkable property of this series is
that for certain values of $s$ and all $n \ge 13$
one can add points to the integer affine
span of $P(n,s)$ in such a way that $P(n,s)$ remains a perfect Delaunay
polytope in the new lattice. Thus, we have constructed an inhomogeneous
analog of the remarkable relationship between $\sfA_9$ and  $\sfA_9^2$.
\end{abstract}

\section{Introduction}

Given a $n$-dimensional lattice $L$, a polytope $D$ is called a
{\em Delaunay polytope} if the set of its vertices is $S\cap L$ with $S$
being a sphere containing no lattice points in its interior.
If $(v_1, \dots, v_n)$ is a basis of $L$ then the Gram matrix $Q=(\langle v_i, v_j\rangle)_{1\leq i,j\leq n}$ characterizes $L$ up to isometry.
It has long been observed that for computations it is preferable
to work with Gram matrices instead of lattices.
Then one defines $S^{n}_{>0}$ the cone of positive definite
$n\times n$-symmetric matrices, identifies the quadratic forms with symmetric
matrices and defines $A[X]=X^{t} AX$ for a column vector $X$ and a symmetric
matrix $A$.

Voronoi \cite{Vo} remarked that if $D$ is a polytope with coordinates in $\ZZ^n$ then the condition that $D$ is a Delaunay polytope is expressed by linear equalities and inequalities on the coefficients of the Gram matrix.
That is if one defines
\begin{equation*}
\SC(D)=\left\{\begin{array}{c}
Q\in S^n_{>0}\quad :\quad \exists c\in\RR^n, r > 0\mbox{~such~that~}\\
Q[v-c]=r\mbox{~for~}v\in\vertt D\\
\mbox{~and~}Q[v-c]>r \mbox{~for~}v\in\ZZ^n-\vertt D
\end{array}\right\}
\end{equation*}
then $\SC(D)$ (called {\em Baranovskii cone} in \cite{AchillBook}) is a polyhedral cone. The dimension of $\SC(D)$ is called the {\em rank} of $D$.
$D$ is called {\em perfect} if it is of rank $1$ (see \cite{Erdahl92} and \cite{combinatorica93} for more details).

The only perfect Delaunay polytope of dimension $n\leq 6$ are the interval $[0,1]$ and Schl\"afli polytope $2_{21}$, which are Delaunay polytopes of the root lattices $\mathsf{A}_1$ and $\mathsf{E}_6$ (see \cite{DD04}).
Several infinite series of perfect Delaunay polytopes were built in \cite{ErdahlInfinite}, \cite{DutourInfinite} and \cite{GrishInfinite}.
Some, conjectured to be complete, lists are given in \cite{IntegerPaper} for dimension $7$ and $8$.
A polytope $P$ is called {\em centrally symmetric} if there exist a point $c$, called {\em center}, such that for any vertex $v\in P$ we have $2c-v\in P$.
In this paper for every $4s\leq n+1$, we build a Delaunay polytope $P(n,s)$ such that:
\begin{enumerate}
\item[(i)] $P(n,s)$ has dimension $n$, is centrally symmetric and has $2{n+1 \choose s}$ vertices.
\item[(ii)] $P(n,s)$ is perfect for $s\geq 2$.
\end{enumerate}
Given a Delaunay polytope $P$ in a lattice $L$,
we denote by $L(P)$ the set of lattice points
that can be expressed as integral sum of vertices of $P$.
$P$ is {\em generating} if $L(P)$ coincides with $L$.

All perfect Delaunay polytopes known so far were generating and the main interest of $P(n,s)$ is that if 
\begin{equation*}
6s < \left\{\begin{array}{rl}
n+1 & \mbox{if $n$ is odd},\\
n & \mbox{if $n$ is even},
\end{array}\right.
\end{equation*}
then there exists a lattice $L'$
such that $P(n,s)$ is a Delaunay polytope in $L'$ and $L(P) \not= L'$.

The polytope $P(7,2)$ is the Gosset polytope $3_{21}$, which is a Delaunay polytope of the root lattice $\mathsf{E}_7$ and 
$P(8,2)$ is the Delaunay polytope $D^8_2$ of \cite{IntegerPaper}.
The infinite series $P(n,s)$ were found by looking at $D^8_2$
and the lattice $L'$ was found by an exhaustive search using the
computer package \cite{polyhedral}.

\section{The lattice $\mathsf{A}_n$}

The lattice $\mathsf{A}_n$ is defined as
\begin{equation*}
\mathsf{A}_n=\left\{x=(x_0,\dots,x_{n})\in \ZZ^{n+1}\quad :\quad \sum_{i=0}^n x_i=0\right\}.
\end{equation*}
$\mathsf{A}_n$ is an $n$-dimensional lattice, but
best seen as embedded into $\RR^{n+1}$ with the standard Euclidean metric $\sum_{i=0}^{n} x_i^2$. 
Define $(e_i)_{1\leq i\leq n+1}$ the standard basis of $\RR^{n+1}$.

It is often useful to think of $\mathsf{A}_n$ as a point lattice.
More formally, define

\begin{equation*}
V_{n,s}=\left\{x=(x_0,\dots,x_{n})\in \ZZ^{n+1}\quad :\quad \sum_{i=0}^{n} x_i=s\right\}.
\end{equation*}
Then the difference set $V_{n,s}-V_{n,s}$ is the lattice $\sfA_n$. Let
\begin{equation*}
J(n+1,s)=\conv\left\{x\in \{0,1\}^{n+1}\quad:\quad\sum_{i=0}^{n} x_i=s\right\}.
\end{equation*}
It is easily seen that $J(n+1,s)$ is a lattice polytope in the point lattice $V_{n,s}$.
 Since $V_{n,s}-V_{n,s}=\sfA_n$, we know that  $\sfA_n$ contains lattice polytopes isometric to $J(n+1,s)$.

For $\alpha_0, \dots, \alpha_n\in \RR$, we define 
\begin{equation*}
q_{\alpha_0,\dots,\alpha_n}(x)=\sum_{i=0}^{n} \alpha_i x_i^2
\end{equation*}
and denote by ${\mathcal QP}$ the cone of all $q_{\alpha_0,\dots,\alpha_n}$
with $\alpha_i>0$.
Clearly the polytopes $J(n+1,s)$ are Delaunay polytopes 
of $\mathsf{A}_n$ for the scalar product induced by $q\in {\mathcal QP}$.

The following theorem is a reformulation of Proposition 8 of \cite{Grish}.
\begin{theorem}\label{ANexample}

(i) The lattice $\mathsf{A}_n$ has $n$ translation classes of Delaunay
polytopes. These classes are represented by polytopes $J(n+1,s)$ for $1\leq
s\leq n$.

(ii) The scalar products on $\mathsf{A}_n$ having the polytopes $J(n+1,s)$
as Delaunay polytopes are the ones induced by some $q\in {\mathcal QP}$.

\end{theorem}
According to the terminology of \cite{Grish} this means that the non-rigidity
degree of $\mathsf{A}_n$ is $n+1$.
Note that the forms $x_0^2,\dots,x_n^2$ remain independent when restricted
to $\mathsf{A}_{n}$.
One classic example is the Delaunay tessellation of $\mathsf{A}_3$: It is formed by the regular simplex $J(4,1)$, its antipodal $J(4,3)$ and the regular octahedron $J(4,2)$.

Clearly, the rank of the polytopes $J(n+1,1)$ and $J(n+1, n)$ is
$\frac{n(n+1)}{2}$ since
those polytopes are $n$-dimensional simplices.

\begin{theorem}\label{FactoidsJNS}
Let $n,s \in \N$ and $2\leq s\leq n-1$.

(i) The rank of $J(n+1,s)$ is $n+1$ and every scalar product on $\mathsf{A}_n$
having $J(n+1,s)$ as Delaunay is induced by some $q\in {\mathcal QP}$.

(ii) The center $c_{\alpha_0,\dots, \alpha_n}$ of the empty ellipsoid
around $J(n+1,s)$ with respect to the quadratic form
$q_{\alpha_0, \dots, \alpha_n}$ is given by
\begin{equation*}
c_{\alpha_0, \dots,\alpha_n}=\left(\frac{1}{2}+\frac{C}{\alpha_0}, \dots, \frac{1}{2}+\frac{C}{\alpha_n}\right)\mbox{~with~}
C=\frac{s-\frac{n+1}{2}}{\sum_{i=0}^{n} \frac{1}{\alpha_i}}.
\end{equation*}
\end{theorem}
\proof For $i=1,\dots,n$ define $v_i=e_i - e_0$.
The norm of a vector $x=\sum_{i=0}^{n} x_i e_i\in \mathsf{A}_n$ with respect to $q_{\alpha_0, \dots,\alpha_n}$ is
\begin{equation*}
\begin{array}{rcl}
q_{\alpha_0,\dots,\alpha_n}(x)
&=& q_{\alpha_0,\dots,\alpha_n}(-(\sum_{i=1}^{n} x_i) e_0+\sum_{i=1}^{n} x_ie_i)\\
&=& \alpha_0(\sum_{i=1}^n x_i)^2+\sum_{i=1}^{n} \alpha_i x_i^2\\
&=& X^{t} A_{\alpha_{0},\dots,\alpha_n} X
\end{array}
\end{equation*}
where $X=(x_1,\dots, x_n)^{t}$, and
\begin{equation*}
A_{\alpha_0,\dots,\alpha_n}=\left(\begin{array}{cccc}
\alpha_0+\alpha_1 &\alpha_0 &\dots   &\alpha_0\\
\alpha_0          &\alpha_0+\alpha_2 &\ddots  &\vdots\\
\vdots            &\ddots            &\ddots  &\alpha_0\\
\alpha_0          &\dots             &\alpha_0&\alpha_0+\alpha_n
\end{array}\right).
\end{equation*}
Expressed in terms of the basis $(v_i)_{1\leq i\leq n}$ the polytope $J(n+1,s)$ is written as
\begin{equation*}
J'(n+1,s)=\conv\left\{(x_1,\dots,x_n)\in \{0,1\}^n\mbox{~with~}s-\sum_{i=1}^nx_i\in \{0,1\}\right\}.
\end{equation*}
Theorem \ref{ANexample}, (ii) then implies that if $\alpha_i>0$, then $A_{\alpha_0,\dots,\alpha_n}\in \SC(J'(n+1,s))$.
Let us now take $A=(a_{i,j})_{1\leq i,j\leq n}\in \SC(J'(n+1,s))$.

Select a three element subset $S=\{s_1,s_2,s_3\}$ of $ \{1,\dots, n\}$ and a vector $v\in
\{0,1\}^{n}$. Consider the polytope 
\begin{equation*}
J_{S, v}=\conv\{w\in \vertt J'(n+1,s)\quad :\quad w_i=v_i\mbox{~for~}i\notin S\}.
\end{equation*}
If one chooses $v$ such that $\sum_{i\notin S} v_i=s-2$, then $J_{S, v}$ is affinely equivalent to the octahedron $J(4,2)$.
The quadratic form $q(x)=X^{t} AX$ induces a quadratic form $q_S$ on the affine space spanned by $J_{S,v}$ with $q_{S}(Y)=Y^{t} A_{S} Y$, $Y=(x_{s_1}, x_{s_2}, x_{s_{3}})^t$ and $A_{S}=(a_{i,j})_{i,j\in S}$.

The rank of the octahedron $J(4,2)$ is equal to $4$ as proved on page 232 of \cite{DL}.
The quadratic form $A_{\alpha_0, \alpha_{1}, \alpha_{2}, \alpha_{3}}$ with $\alpha_i>0$ has $4$
independent coefficients and belongs to $\SC(J_{S,v})$ thus we get
$A_S=A_{\alpha_0,\dots,\alpha_3}$ for some coefficients $\alpha_i$.
This implies $a_{i,j}=C_{S}$ for $i\not=j \in S$ with the
constant $C_{S}$ a priori depending on $S$.
If one interprets the value $a_{i,j}$ as colors of an edge between vertices
$i$ and $j$ then we get that all triangles of the complete graph on $n$
vertices are monochromatic. This is possible only if there is only one edge
color. 
So, $a_{i,j}=C$ for $i\not= j$ and one can write $A=A_{\alpha_0, \dots, \alpha_n}$ with $\alpha_i\in \RR$.

Let us find the circumcenter of the empty sphere around $J(n+1,s)$.
The point $h_{n+1}=((\frac{1}{2})^{n+1})$ is at equal distance from all points of
$J(n+1,s)$. However, it does not belong to $V_{n,s}$.
To find the circumcenter $c$ of $J(n+1,s)$, we take the orthogonal
projection of $h_{n+1}$ on the hyperplane $\sum_{i=0}^{n} x_i=s$ for the quadratic
form $q_{\alpha_0, \dots, \alpha_n}$. Easy computations give (ii).

Let us prove $\alpha_i>0$.
It is well known that the facets of $J(n+1,s)$ are determined by the inequalities $x_i\geq 0$ and $x_i\leq 1$. 
It is also easy to see that the Delaunay polytopes adjacent to the facets $x_0\geq 0$ and $x_0\leq 1$ are
\begin{equation*}
\begin{array}{rcl}
J_0^{-}&=&\{x\in \{-1, 0\}\times \{0,1\}^n\quad :\quad \sum_{i=0}^{n} x_i=s\},\\
\mbox{and~}J_0^{+}&=&\{x\in \{1, 2\}\times \{0,1\}^n\quad :\quad \sum_{i=0}^{n} x_i=s\}.
\end{array}
\end{equation*}
The polytopes $J_0^{-}$, $J_0^{+}$ are equivalent under translation to $J(n+1,s+1)$ and $J(n+1,s-1)$.

The square distance of $h_{n+1}$ to the vertices of $J(n+1,s)$ is
$d=\sum_{i=0}^n \frac{\alpha_i}{4}$ and the square distance of $h_{n+1}$
to the vertices of $J_0^{-}$, $J_0^{+}$ not in $J(n+1,s)$
is $d'=\alpha_0 \frac{9}{4}+\sum_{i=1}^n \alpha_{i}\frac{1}{4}$.
The conditions defining $\SC(J'(n+1,s))$ imply $d'>d$ hence
$\alpha_0>0$ and by symmetry $\alpha_i>0$.
So, the conditions for $J(n+1,s)$ to be a Delaunay polytope imply that
$A=A_{\alpha_0,\dots, \alpha_n}$ with $\alpha_i>0$.
But according to Theorem \ref{ANexample} those conditions are sufficient
for the stronger condition of preserving all the Delaunay polytopes of
${\mathsf A}_n$ so they are clearly sufficient for just $J(n+1,s)$. \qed

\section{The polytopes $P(n,s)$}\label{CompMethods}
We denote an $(n+1)$-vector whose first $a$ coordinates are $A$ and the remaining $n+1-a$ coordinates $B$ by $(A^a;B^{n+1-a})$. Similar convention is used for vectors with three distinct coordinates, e.g.  $(A^a;B^b;C^{n+1-a-b})$.

\begin{definition}
Take $n,s\in \ZZ$ with $s\geq 1$ and $4s\leq n+1$.

(i) Set $v_{n,s}=\left(\left(\frac{1}{4}\right)^{4s}; 0^{n+1-4s}\right)$.
The polytope $P(n,s)$ is defined as 
\begin{equation*}
P(n,s)=\conv\left\{v, 2v_{n,s}-v \mbox{~for~}v\in \vertt J(n+1,s)\right\}.
\end{equation*}

(ii) Define $t_{n,s}=\left(\left(\frac{1}{2}\right)^{2s};  \left(\frac{-1}{2}\right)^{2s};  0^{n+1-4s}\right)$ and 
\begin{equation*}
V^2_{n,s}=\left\{v, t_{n,s}+ v \mbox{~for~} v\in V_{n,s}\right\}.
\end{equation*}

\end{definition}

\begin{theorem}\label{InfiniteSeriesBirth}
Take $n,s\in \ZZ$ with $s\geq 2$ and $4s\leq n+1$.

(i) $V^{2}_{n,s}$ is a lattice and $P(n,s)$ affinely generates it.

(ii) The polytope $P(n,s)$ is perfect with the unique, up to positive multiple, positive definite quadratic form being 
\begin{equation*}
q_{n,s}(x)=2\sum_{i=0}^{4s-1} x_i^2+\sum_{i=4s}^{n} x_i^2.
\end{equation*}
The center of the circumscribed ellipsoid is $v_{n,s}$ and the squared radius is $\frac{3s}{2}$.

\end{theorem}
\proof We have $2t_{n,s}\in V_{n,s}$ so $V^2_{n,s}$ is a lattice. $P(n,s)$ generates it since $J(n+1,s)$ generates $\mathsf{A}_n$.
By its definition,  $P(n,s)$ is centrally symmetric of center $v_{n,s}$.
It is well known and easy to prove that if a Delaunay polytope is centrally
symmetric then the center $c'$ of its empty sphere coincides with the center
$c$ of the antisymmetry operation $v\mapsto 2c-v$.
So, we should have $v_{n,s}=c_{\alpha_0,\dots, \alpha_n}$ with
\begin{equation*}
c_{\alpha_0, \dots,\alpha_n}=\left(\frac{1}{2}+\frac{C}{\alpha_0}, \dots, \frac{1}{2}+\frac{C}{\alpha_n}\right).
\end{equation*}
Thus:
\begin{itemize}
\item For $0\leq i\leq 4s-1$, we have $c_i=\frac{1}{4}$.
This implies $\alpha_i=-4C$.
\item For $4s\leq i\leq n$, we have $c_i=0$. 
This implies $\alpha_i=-2C$.
\end{itemize}
Summarizing we get $q=-2C q_{n,s}$ and thus that $P(n,s)$ is perfect.
The proof of the Delaunay property follows from the fact that the coefficient in front of $x_i^2$ are strictly positive for $0 \le i \le n$ and property (i) of Theorem \ref{FactoidsJNS}. \qed

\section{The lattice $V^4_{n,s}$}

Define the vector $w_{n,s}$ by
\begin{equation*}
w_{n,s}=\left\{\begin{array}{rl}
\left(\left(\frac{1}{4}\right)^{2s} , \left(\frac{-1}{4}\right)^{2s}, \left(\frac{1}{2}\right)^{n+1-4s}\right)-\frac{n+1-4s}{2} e_1 & {\rm if}\, n\, {\rm is~odd}, \\[0.8mm]
\left(\left(\frac{1}{4}\right)^{2s} , \left(\frac{-1}{4}\right)^{2s}, 0, \left(\frac{1}{2}\right)^{n-4s}\right)-\frac{n-4s}{2} e_1 & {\rm if}\, n\, {\rm is~even}.
\end{array}\right.
\end{equation*}
Then define
\begin{equation*}
V^4_{n,s}=V^2_{n,s}\cup w_{n,s}+V^2_{n,s}.
\end{equation*}
Clearly $V^4_{n,s}$ is a lattice that contains $V^2_{n,s}$ as an index $2$ sublattice. We want to prove that $P(n,s)$ remains a Delaunay polytope in $V^4_{n,s}$ for some values of $n$ and $s$.

\begin{theorem}\label{GeneratingTheor}
The polytope $P(n,s)$ is a Delaunay polytope of $V^4_{n,s}$ if
\begin{equation*}
6s < \left\{\begin{array}{rl}
n+1 & \mbox{if $n$ is odd},\\
n & \mbox{if $n$ is even}.
\end{array}\right.
\end{equation*}

\end{theorem}
\proof We need to solve the closest vector problem for the
lattice $V^4_{n,s}$ and the point $v_{n,s}$.
For $V^2_{n,s}$ this is solved by Theorem \ref{InfiniteSeriesBirth}.
Thus we need to find the closest vectors in $w_{n,s}+V^{2}_{n,s}$ to $v_{n,s}$.
This is equivalent to finding the closest vectors in $V_{n,s}$ to $v_{n,s} - w_{n,s}$ and to $v_{n,s} - w_{n,s} - t_{n,s}$.
We have if $n$ is odd:
\begin{equation*}
\begin{array}{rcl}
v_{n,s} - w_{n,s}
&=& \left(0^{2s}; \left(\frac{1}{2}\right)^{2s}; \left(-\frac{1}{2}\right)^{n+1-4s}\right)+ \frac{n+1-4s}{2} e_1,\\[0.5mm]
v_{n,s} - w_{n,s} - t_{n,s}
&=& \left(\left(-\frac{1}{2}\right)^{2s}; 1^{2s}; \left(-\frac{1}{2}\right)^{n+1-4s}\right)+ \frac{n+1-4s}{2} e_1,
\end{array}
\end{equation*}
and if $n$ is even:
\begin{equation*}
\begin{array}{rcl}
v_{n,s} - w_{n,s}
&=& \left(0^{2s}; \left(\frac{1}{2}\right)^{2s}; 0; \left(-\frac{1}{2}\right)^{n-4s}\right) + \frac{n-4s}{2} e_1,\\[0.5mm]
v_{n,s} - w_{n,s} - t_{n,s}
&=& \left(\left(-\frac{1}{2}\right)^{2s}; 1^{2s}; 0; \left(-\frac{1}{2}\right)^{n-4s}\right) + \frac{n-4s}{2} e_1.
\end{array}
\end{equation*}
All the vectors occurring have coordinates belonging to $\ZZ$ or
$\ZZ+\frac{1}{2}$. Since the coordinates of elements of $V_{n,s}$ are
integral and $q_{n,s}$ has non-zero coefficients only for $x_i^2$
this gives for $v\in V_{n,s}$ the following lower bounds if $n$ is odd:
\begin{equation*}
\begin{array}{rcl}
q_{n,s}(v_{n,s} - w_{n,s} -v)
&\geq & 2\times 2s\times \frac{1}{4} + (n+1-4s)\frac{1}{4}=\frac{n+1}{4},\\[1mm]
q_{n,s}(v_{n,s} - w_{n,s} -t_{n,s} -v)
&\geq & 2\times 2s\times \frac{1}{4} + (n+1-4s)\frac{1}{4}=\frac{n+1}{4},
\end{array}
\end{equation*}
and if $n$ is even:
\begin{equation*}
\begin{array}{rcl}
q_{n,s}(v_{n,s} - w_{n,s} - v)
&\geq & 2\times 2s\times \frac{1}{4} + (n-4s)\frac{1}{4} =\frac{n}{4},\\[1mm]
q_{n,s}(v_{n,s} - w_{n,s} - t_{n,s} - v)
&\geq & 2\times 2s\times \frac{1}{4} + (n-4s)\frac{1}{4} =\frac{n}{4}.
\end{array}
\end{equation*}
So, if $n$ and $s$ satisfy the condition of the theorem then the closest
points in $w_{n,s} + V^2_{n,s}$ are at a square distance greater than
$\frac{3s}{2}$. But $\frac{3s}{2}$ is the square radius of the
circumscribing sphere thus proving that $P(n,s)$ is a Delaunay
polytope in $V^{4}_{n,s}$. \qed

The above theorem gives example of Delaunay polytopes, which are perfect but not generating, the first example of which is $P(13,2)$.

\end{document}